\documentclass[oribibl]{llncs}

\def\Z{\bbbz} %Integers
\def\Q{\bbbq} %Quotients
\def\F{\bbbf} %Finite Field(nonprescribe)
\def\Fq{\bbbf_q} %FiniteField
\def\lg{{\mathrm{lg}\,}}
\def\lc{{\mathrm{lc}\,}}

\usepackage{enumerate}
\spnewtheorem{algorithm}[theorem]{Algorithm}{\bfseries}{\rmfamily}

\begin{document} 
\title{More constructing pairing-friendly elliptic curves for cryptography}
\author{Tanaka Satoru \and Nakamula Ken}
\institute{Department of Mathematics and Information Sciences, Tokyo Metropolitan University, \\
1-1 Minami Osawa, Hachioji-shi, Tokyo, 192-0397 Japan\\
\email{satoru@tnt.math.metro-u.ac.jp, nakamula@tnt.math.metro-u.ac.jp}
}

%\date{}

\maketitle

\begin{abstract}{
  The problem of constructing elliptic curves suitable for pairing applications has received a lot of attention. To solve this, we propose a variant algorithm of a known method by Brezing and Weng. We produce new families of parameters using our algorithm for pairing-friendly elliptic curves of embedding degree 8, and we actually compute some explicit curves as numerical examples.
%, for instance, to compute the Ate Pairing \cite{HSV06}. %in Section \ref{sec:Res_Comp_Result}.
}
\end{abstract}

\section{Introduction}
\label{sec:Intro}

  Researches on pairing-based cryptographic schemes have received interest over the past few years. Recently many new and novel protocols have been proposed as in \cite{SOK00,BF03,BLS01,Jou04}. A randomly chosen elliptic curve, however, rarely has a subgroup of large prime order, therefore construction of ``pairing-friendly'' elliptic curves is one of the important problems for cryptography \cite{BK98}.

  Let $E$ be an elliptic curve defined over a finite field $\Fq$, and $r$ be the largest prime dividing $\#E(\Fq)=q+1-t$, the order of the group of $\Fq$-rational points of $E$ with the Frobenius trace $t$. We define the {\it embedding degree}
%, which is the degree $k$ of the extension field required by the pairing,
 by the smallest positive integer $k$ such that $r$ divides $q^k-1$. 
%Besides, we use ``CM discriminant'' $D$ by CM method\cite{AM93}. 
%The parameters required for the construction of pairing-friendly elliptic curves are $t,r,q,k$ and $D$ and we will use the so-called CM method introduced in \cite{AM93} to solve for the parameters.
The parameters required to determine pairing-friendly elliptic curves are $t,r,q,k$ and the CM discriminant $D$ for the CM method introduced in \cite{AM93} to construct elliptic curves.

  In this paper,  we study the problem of computing suitable parameters $t, r, q$ from given parameters $k, D$. We employ the method proposed in \cite{FST06,BW05} which generates a family of pairing-friendly curves by considering $t, r, q$ as polynomial $t(x), r(x), q(x)$ of a new parameter $x$.
% We modified the original method by changing the input to finite subset of $k$th cyclotomic field $\Q(\zeta_k)$.
 We restrict the embedding degree to $k=8$ and the CM discriminant to $D=1$. The key point is how to choose a good $r(x)$. Instead of taking $r(x)$ to be the $\ell$th cyclotomic polynomial $\Phi_\ell(x)$ with a multiple $\ell$ of $k$, we modified the original method by starting from a finite subset of the $k$-th cyclotomic field $\Q(\zeta_k)$ with a primitive $k$th root $\zeta_k$ of unity so that $r(x)$ can be systematically computable.
 As a result, we came up with new families of pairing-friendly curves which are given explicitly in Table~\ref{tab:family_cubicdeg8} and Theorem~\ref{thm:cubic8_t82} of Section~\ref{sec:Pair_Res}. We also give, for the first time in this case so far as we know, explicit numerical results as in Examples 1--3. %We are now able to systematically construct pairing-friendly elliptic curves, so we expect that pairing-based cryptosystems will be constructed by selecting pairing-friendly curves using our algorithm. 
%Also, we will put new pairing-based cryptographic scheme to practical use.

%  $B$3$NO@J8$G07$&%Z%"%j%s%0$N7W;;$KE,$7$?BJ1_6J@~$N9=@.LdBj$O!$M?$($i$l$?%Q%i%a!<%?(B$k,D$$B$KBP$7$F!$(BCM$BK!$rMQ$$$F6J@~$r9=@.$9$k$N$KI,MW$J%Q%i%a!<%?(B$t,r,q$$B$r7W;;$9$k$3$H$G$"$k!%2f!9$O(B\cite{FST06,BW05}$B$GDs0F$5$l$?!$(B$t,r,q$$B$r?7$7$$%Q%i%a!<%?(B$x$$B$NB?9`<0$NAH(B$t(x),r(x),q(x)$$B$H9M$(!$%Z%"%j%s%0$KE,$7$?6J@~$N%U%!%_%j$r@8@.$9$kJ}K!$r:NMQ$9$k!%:#2s!$Kd$a9~$_<!?t$r(B$k=8$$B!$(BCM$BH=JL<0$r(B$D=1$$B$K8GDj$7!$?7$?$K(B\S~\ref{sec:Pair_Res}$B$N(B{\bf $B%"%k%4%j%:%`(B~\ref{algo:cubic8}}$B$K$h$C$F(B$t(x),r(x),q(x)$$B$r7W;;$7$?!%$=$N7k2L!$(B\S~\ref{sec:Res_Comp_Result}$B$N(B{\bf $BDjM}(B~\ref{thm:cubic8_t82}}$B$K<($5$l$k$h$&$J?7$7$$%Z%"%j%s%0$KE,$7$?6J@~$N%U%!%_%j$N9=C[$K@.8y$7$?!%K\O@J8$GDs0F$9$k<jK!$K$h$j!$4{B8$N<jK!$G$OC5:w$,:$Fq$G$"$C$?6J@~$N%U%!%_%j$N9=@.$rMF0W$K9T$($k$h$&$K$J$k!%$3$N7k2L$H$7$F!$%Z%"%j%s%0$KBP$7$FM%$l$?6J@~$r9=@.$9$k$3$H$,2DG=$K$J$k$@$1$G$J$/!$$3$NJ}<0$GF@$i$l$?6J@~$rMQ$$$?%Z%"%j%s%0$K$h$k0E9f7O$N9=@.$,9T$($k$h$&$K$J$k$H4|BT$5$l$k!%(B

  This paper is organized as follows. Section~\ref{sec:Pair_Frame} gives a brief mathematical definition of curves suitable for pairing-based cryptography and the method of construction we used to generate our curves. In Section~\ref{sec:Pair_Res} we give our algorithm to construct curves. Section~\ref{sec:Res_Comp_Result} gives numerical examples of curves that we generate using our parameters. Finally, we will discuss the conclusions that we will draw from our approach in Section~\ref{sec:Cons}.

\section{Our framework of pairing-friendly curves}
\label{sec:Pair_Frame}

  A survey on the construction of pairing-friendly elliptic curves is given by Freeman et al. \cite{FST06}. We introduce several essential definitions from that paper to explain our algorithm.
% We use Freeman's denotion without notice and use denotion $\lg := \log_2$. 
We will use the same notation there without notice. Let $\lg$ mean the base 2 logarithm in the following.

\subsection{Families of curves for pairing}
\label{sec:Pair_Family}

  At first, we give the definition of pairing-friendly elliptic curves used in cryptography.
\begin{definition}[\protect{\cite[Definition 2.3]{FST06}}]
\label{def:PairingFriendly}\rm
  Suppose $E$ is an elliptic curve defined over $\Fq$. We say that $E$ is {\it pairing-friendly} if $E$ satisfies the following conditions:
  \begin{enumerate}
    \item \label{cond:Pair_scale} there is a prime $r \geq \sqrt{q}$ such that $r \mid \#E(\Fq)$.
    \item \label{cond:Pair_embdeg} the embedding degree of $E$ with respect to $r$ is less than $(\lg r)/8$.
  \end{enumerate}
\end{definition}
 For cryptographic applications of pairings, basically we desire enough security depending on the {\it elliptic curve discrete logarithm problem} ({\it ECDLP}). In fact, by this definition, suitable sizes of $r, q^k$ seem to avoid any known attack for the ECDLP today \cite{FST06}.

  Next, we explain how to construct pairing-friendly curves. The parameter $q$ has to be a prime power. If a family of pairing-friendly curves represented by $t(x), r(x)$ and $q(x)$ is 
given, we anticipate that 
% not only several $x$ satisfies that the $q(x)$ is prime power but 
%there are infinitely many $x$ satisfying the condition that
 $q(x)$ is a prime power for infinitely many $x$. Freeman et al. gave a definition with a familiar conjecture as follows \cite[Section~2]{FST06}.

\begin{definition}
\label{def:repr_primes}\rm
  Let $f(x)$ be a polynomial with rational coefficients. We say $f$ {\it represents primes} if the following conditions are satisfied.
  \begin{enumerate}
    \item \label{cond:rep_irred} $f(x)$ is non-constant and irreducible.
    \item \label{cond:rep_led} $f(x)$ has positive leading coefficient. 
    \item \label{cond:rep_zahl} $f(x) \in \Z$ for some $x \in \Z$.
    \item \label{cond:rep_coprime} $\gcd(\{f(x) \mid x,f(x) \in \Z\})=1$.
  \end{enumerate}
\end{definition}

  We use this definition to define families of pairing-friendly curves.
\begin{definition}[\protect{\cite[Definition 2.6]{FST06}}]
\label{def:pair_friendly_family}\rm
  For a given positive integer $k$ and positive square-free integer $D$, the triple {\it $(t,r,q)$ represents a family of elliptic curves with embedding degree $k$ and discriminant $D$} if the following conditions are satisfied:
  \begin{enumerate}
    \item \label{cond:family_FF} $q(x)=p(x)^d \,\,(d \ge 1)$ and $p(x)$ that represents primes.
    \item \label{cond:family_ECDLP} $r(x)=c \cdot {\tilde{r}}(x) \,\,(c \in \Z_{\ge 1})$ and $\tilde{r}(x)$ that represents primes.
    \item \label{cond:family_subgrp} $r(x) \mid q(x)+1-t(x)$.
    \item \label{cond:family_embdeg_cyclo}  $r(x) \mid \Phi_k(t(x)-1)$, where $\Phi_k$ is the $k$th cyclotomic polynomial.
    \item \label{cond:family_CMeq} The CM equation $4q(x)-t(x)^2 = Dy^2\,\,(D \in \Z_{>0})$ has infinitely many integer solutions $(x,y)$.
  \end{enumerate}
\end{definition}

  For a family $(t(x),r(x),q(x))$, if the CM equation in (\ref{cond:family_CMeq}) has a suitable set of integer solutions $(x_0,y_0)$ with both of $p(x_0)$ and $\tilde{r}(x_0)$ are primes, then we are able to construct curves $E$ over $\F_{q(x_0)}$ where $E(\F_{q(x_0)})$ has a subgroup of order $\tilde{r}(x_0)$ and embedding degree $k$ with respect to $\tilde{r}(x_0)$ by using the CM method in \cite{AM93}. %After check Definition~\ref{def:PairingFriendly}, we will get pairing-friendly elliptic curves.

  We therefore define a parameter $\rho$ that represent how close to the ideal curve that is $\#E(\Fq)$ is prime as follows.
\begin{definition}[\protect{\cite[Definition 2.7]{FST06}}]
\label{def:rhovalueofFamily}\rm
\begin{enumerate}
  \item Let $E/\Fq$ be an elliptic curve, and suppose $E$ has a subgroup of order $r$. The {\it $\rho$-value of $E$ (with respect to $r$)} is
    \[
      \rho(E) = \frac{\log q}{\log r} .
    \]
  \item Let $t(x),r(x),q(x)\,\in \Q[x]$, and suppose $(t, r, q)$ represents a family of elliptic curves with embedding degree $k$. The {\it $\rho$-value of the family represented by $(t,r,q)$} is
    \[
      \rho(t,r,q) = \lim_{x \to \infty}\frac{\log q(x)}{\log r(x)} = \frac{\deg q(x)}{\deg r(x)} .
    \]
\end{enumerate}
\end{definition}
  By Definition~\ref{def:PairingFriendly}, a pairing-friendly curve $E$ has $\rho(E) \le 2$. The smaller the $\rho$-value, the faster the computation of points on elliptic curve (See \cite[Section 1.1]{FST06}). On the other hand, the Hasse bound implies that $\rho(t,r,q)$ is always at least $1$.
%so $r \approx q$ -- to find an elliptic curve $E$ satisfies that $\rho(E) \approx 1$ more efficiently are important problem for cryptography.
Finding parameters efficiently with the same bit size of $r$ and $q$, hence $\rho(E)$ is close to $1$, is one of the important problems for cryptography. 

\subsection{Original method}
\label{sec:Res_BWMethod}

%  For example of Definition \ref{def:pair_friendly_family}, 
  In this section, we briefly explain the construction of curves satisfying the condition of Definition \ref{def:pair_friendly_family} proposed by Brezing and Weng \cite{BW05}\cite[Section 6.1]{FST06}.

\begin{theorem}
\label{thm:ExtCP}
  Fix a positive integer $k$ and a positive square-free integer $D$. Execute the following steps:
  \begin{enumerate}[{\rm Step 1.}]
    \item \label{algstep:BW_selfld}Choose a number field $K$ containing a primitive $k$th root of unity $\zeta_k$ and $\sqrt{-D}$.
    \item \label{algstep:BW_selr}Find an irreducible polynomial $r(x) \in \Z[x]$ such that $\Q[x]/(r(x)) \cong K$.
    \item \label{algstep:BW_selt}Let $t(x) \in \Q[x]$ be a polynomial mapping to $\zeta_k +1 \in K$.
    \item \label{algstep:BW_sely}Let $y(x) \in \Q[x]$ be a polynomial mapping to $(\zeta_k -1)/\sqrt{-D} \in K$. (So, if we discover a polynomial $s(x) \in \Q[x]$ mapping to $\sqrt{-D} \in K$, then $y(x) \equiv (2-t(x))s(x)/D \pmod{r(x)}$.)
    \item \label{algstep:BW_selq}Let $q(x)=({t(x)}^2+D{y(x)}^2)/4$.
  \end{enumerate}
  If both of $r(x)$ and $q(x)$ represent primes, then the triple $(t,r,q)$ represents a family of curves with embedding degree $k$ and CM discriminant $D$.
\end{theorem}

  The $\rho$-value for this family is
\[
  \rho(t,r,q) = \frac{2\max\{\deg t(x), \deg y(x)\}}{\deg r(x)}<2 .
\]
  For more details, refer to \cite[Section 6.1]{FST06}. To find a family of pairing-friendly elliptic curves efficiently, we have to choose a good $r(x)$ satisfying $\zeta_k,\sqrt{-D} \in K$. The idea by Brezing, Weng and also Freeman et al. is as follows. Choose an integer multiple $\ell$ of $k$ so that $\sqrt{-D} \in K = \Q(\zeta_\ell)$. Then let $r(x)=\Phi_\ell(x)$. %Since $\zeta_k \in K=\Q(\zeta_\ell)$, we can easily check whether $\sqrt{-D} \in K$ .
% with the cyclotomic field theory, so we determined suitable $\ell$.
%Using this idea, we can find a suitable integer multiple $\ell$ of $k$.
They further give some sporadic families \cite[Section 6.2]{FST06}. 
Our idea given explicitly below is to construct such sporadic curves systematically.

\section{Proposed algorithm}
\label{sec:Pair_Res}

\subsection{Factorization of cyclotomic polynomial}
\label{sec:Res_Cube8}

  When we use the original method to construct families, the problem is how to choose polynomials at Step \ref{algstep:BW_selr} and \ref{algstep:BW_selt} in Theorem~\ref{thm:ExtCP}.
If $\Phi_k(u(x))$ has a factorization over $\Q$ for some $u(x) \in \Q[x]$, we let $r(x)$ be one of the irreducible factors. Set $K=\Q[x]/(r(x))$ and we will get $u(x) \mapsto \zeta_k$.
But these factorizations are rare, so this technique to construct families is called ''Sporadic'' families by Freeman \cite[\S6.2]{FST06}. 

  One of the technique to find such $u(x)$ was discussed in Galbraith, Mckee and Valen\c{c}a by proving an important lemma below \cite[Lemma 1]{GMV05}. Baretto and Naehrig \cite{BN06} found a family of embedding degree 12 and $u(x)$ is a quadratic polynomial with $\rho(t,r,q)=1$ using this lemma. It was restricted to quadratic polynomials $u(x)$. In fact, it is effective for general case as is easily seem from the proof there:
\begin{lemma}
\label{lem:subs_cyclo_factor}
  Let $u(x) \in \Q[x]$ and % be a rational polynomial. and $\Phi_k(u(x))$ be the $k$th cyclotomic polynomial with substitute $x$ to $u(x)$.
% The degree of $\Phi_k(u(x))$ is a multiple of $\varphi(k)$ (
$\varphi$ be Euler function.
%).
 Then, the polynomial $\Phi_k(u(x))$ has an irreducible factor of degree $\varphi(k)$ if and only if the equation
\begin{equation}
  u(x)=\zeta_k
  \label{eq:cyclo_factor}
\end{equation}
has a solution in $\Q(\zeta_k)$.
\end{lemma}
%\begin{proof}
%  Obviously proved by \cite{GMV05}, Lemma 1.
%\qed\end{proof}
  We rediscover families of elliptic curves by Freeman\cite[Example 6.18]{FST06} using Lemma~\ref{lem:subs_cyclo_factor} and we try to construct new families of curves. We propose an algorithm for the construction of a family of curves using Lemma~\ref{lem:subs_cyclo_factor} and Theorem~\ref{thm:ExtCP}. The algorithm is as follows.
\begin{algorithm}
\label{algo:Ext_BW} \
\begin{description}
  \item[Input] Positive integers $D, k$ such that $\sqrt{-D} \in \Q(\zeta_k)$ and a finite subset $S \subset \Q(\zeta_k)$.
  \item[Output] Families of elliptic curves with parameters $t(x), r(x), q(x)$.
\end{description}
  \begin{enumerate}[{\rm Step 1.}]
    \item \label{extbw:linearsolve} For each $\omega \in S$, compute $u(x) \in \Q[x]$ such that the equation~(\ref{eq:cyclo_factor}) has a root $x=\omega$. If $u(x)$ does not exist for all $\omega$, then the algorithm fails.% and terminate.
    \item \label{extbw:factor} For each $u(x)$ at Step \ref{extbw:linearsolve}, compute all irreducible factors $r(x)$ of the polynomial $\Phi_k(u(x))$.
    \item \label{extbw:selt} For each pair of $u(x),r(x)$ at Step \ref{extbw:factor}, compute all polynomials $t(x) \in \Q[x]$ such that $\deg t(x)<\deg r(x)$ and $t(x) \equiv u(x)^m+1 \pmod{r(x)}$ for all $m$ with $1 \le m < k, \gcd(m,k)=1$.
    \item \label{extbw:construct} For each pair of $r(x),t(x)$ at Step \ref{extbw:selt}, execute Step \ref{algstep:BW_sely} and Step~\ref{algstep:BW_selq} of Theorem~\ref{thm:ExtCP} to compute $q(x)$.
    \item \label{extbw:check} For each triple $r(x),t(x),q(x)$ at Step~\ref{extbw:construct}, check whether $q(x),r(x)$ represent primes. If $q(x),r(x)$ represent primes, output a family $t(x),r(x),q(x)$.
  \end{enumerate}
\end{algorithm}
%In this algorithm, the condition $\sqrt{-D} \in \Q(\zeta_k)$ do not need between Step~\ref{extbw:linearsolve} and Step~\ref{extbw:selt}.

\subsection{Algorithm refinement with method of indeterminate coefficients}
\label{sec:Res_Ext_BW_Algo}

  Let $\deg u(x) = 3$ as the case $\deg u(x) \le 2$ is studied in \cite{GMV05}. %We apply Algorithm~\ref{algo:Ext_BW} to $u(x)$ and 
Set the embedding degree to $k=8$. In Step~\ref{extbw:linearsolve} of Algorithm~\ref{algo:Ext_BW}, we employ the method of indeterminate coefficients to compute $u(x)$. This technique is also applicable for general $k$.

  Write any rational cubic polynomial $u(x)$ with coefficients $u_0,u_1,u_2,u_3$ as follows:
\begin{equation}
  u(x) = \sum_{i=0}^{3}u_i x^i={u_3}x^3+{u_2}x^2+{u_1}x+u_0 \quad (u_i \in \Q, u_3 \not=0) .
  \label{eq:u_coeff_repr}
\end{equation}
We represent a given value $\omega \in \Q(\zeta_8)$ as follows:
\begin{equation}
  \omega = \sum_{i=0}^{3}a_i {\zeta_8}^i=a_0+{a_1}\zeta_8+{a_2}{\zeta_8}^2+{a_3}{\zeta_8}^3 \quad (a_i \in \Q) .
  \label{eq:omega_zeta8_repr}
\end{equation}
To avoid operation in $\Q(\zeta_8)$, we replace $\zeta_8$ to $x$ to get the following polynomial.
\[
  \omega(x) = \sum_{i=0}^{3}a_i x^i={a_3}x^3+{a_2}x^2+{a_1}x+a_0 .
\]
Next we look at polynomial $u(\omega(x))$. The equation~(\ref{eq:cyclo_factor}) is equivalent to $u(\omega(x))\equiv x \pmod{\Phi_8(x)}$. We take
\[
  v(x) \equiv u(\omega(x)) \pmod{\Phi_8(x)}
\]
be the simplified polynomial of degree not exceeding three with coefficients expressed in terms of $u_i, a_i$. The equation (\ref{eq:cyclo_factor}) is transformed to the polynomial equation
\begin{equation}
  v(x)=x .
  \label{eq:cubic8_reduce}  
\end{equation}
 We can easily show that the coefficients of the left hand side of the equation are all represented as linear combinations of $u_i$. %So we compare coefficients of equation~(\ref{eq:cubic8_reduce}) then we determine $u_i$. We only solve the linear equation system as follows:
More precisely, it is reduced to solve the following system of linear equations to obtain $u_0, u_1, u_2, u_3$.
\begin{equation}
  A\left(
  \begin{array}{c}
   u_0 \\ u_1 \\ u_2 \\ u_3
  \end{array}
  \right)=
  \left(
  \begin{array}{c}
   0 \\ 1 \\ 0 \\ 0
  \end{array}
  \right) ,
  \label{eq:cubic8_line_sys}  
\end{equation}
where
\[
  A=\left(
  \arraycolsep1ex
  \begin{array}{ccrr}
   1 & a_0 & {a_0}^2-{a_2}^2-2{a_1}{a_3} & {a_0}^3-3{a_2}({a_0}{a_2}+{a_1}^2-{a_3}^2)-6{a_0}{a_1}{a_3} \\
   0 & a_1 & 2{a_0}{a_1}-2{a_2}{a_3} & {a_3}^3-3{a_1}({a_1}{a_3}+{a_2}^2-{a_0}^2)-6{a_0}{a_2}{a_3} \\
   0 & a_2 & {a_1}^2-{a_3}^2+2{a_0}{a_2} & -{a_2}^3+3{a_0}({a_0}{a_2}+{a_1}^2-{a_3}^2)-6{a_1}{a_2}{a_3} \\
   0 & a_3 & 2{a_1}{a_2}+2{a_0}{a_3} & {a_1}^3-3{a_3}({a_1}{a_3}+{a_2}^2-{a_0}^2)+6{a_0}{a_1}{a_2}
  \end{array}
  \right) .
  \label{eq:cubic8_line_array}  
\]
%To solve these equations, 
Let $d$ and $n_i$ be as follows:
\begin{equation}
  \begin{array}{rcl}
    d & := & ({a_1}^2+{a_3}^2)(({a_1}-{a_3})^2+2{a_2}^2)(({a_1}+{a_3})^2-2{a_2}^2) , \\
    n_0 & := & -{a_2}(5{a_1}^4{a_3}-5{a_1}^3{a_2}^2+5{a_1}{a_2}^2{a_3}^2-2{a_2}^4{a_3}+3{a_3}^5) , \\
    n_1 & := & {a_1}^5-4{a_1}^3{a_3}^2+9{a_1}^2{a_2}^2{a_3}+{a_1}(2{a_2}^4+3{a_3}^4)+3{a_2}^2{a_3}^3 , \\
    n_2 & := & {a_1}^3{a_2}+3{a_1}{a_2}{a_3}^2-2{a_2}^3{a_3} ,\\
    n_3 & := & {a_3}^3-{a_1}^2{a_3}+2{a_1}{a_2}^2 .
  \end{array}
  \label{eq:cubic8_invariant}
\end{equation}
  If $d$ is nonzero, then we can solve the system  (\ref{eq:cubic8_line_sys}). The solution is
\begin{equation}
  \left\{
    \begin{array}{rcl}
      u_0 & = & -\left({n_3}{a_0}^3+{n_2}{a_0}^2+{n_1}{a_0}-n_0 \right)/d \\
      u_1 & = & \left(3{n_3}{a_0}^2+2{n_2}{a_0}+{n_1} \right)/d \\
      u_2 & = & -\left(3{n_3}{a_0}+{n_2} \right)/d \\
      u_3 & = & -{n_3}/d \\
    \end{array}
  \right. .
  \label{eq:cubic8_coeffs}
\end{equation}
We now have a concrete solution for $k=8$ with $\deg u(x)=3$. %In general $k$ we will get, in fact the solution is seemed to be complicate. We take a Theorem by this result:
Although this method of indeterminate coefficients can be used for any embedding degree $k$, it is not sure wheter the obtained solution is as simple as the ones we discussed. We present the following theorem as a resut of the solutions we computed.
\begin{theorem}
\label{thm:cubic8_lineq}
  For $\omega \in \Q(\zeta_8)$ given by {(\ref{eq:omega_zeta8_repr})},
% by equation{\rm (\ref{eq:omega_zeta8_repr})},
 let $d, n_i$ be as in {(\ref{eq:cubic8_invariant})}. Then, if and only if both $d$ and $n_3$ are nonzero, the equation {(\ref{eq:cyclo_factor})} has a solution $x=\omega$ for a cubic polynomial $u(x) \in \Q[x]$, which is uniquely determined by {(\ref{eq:cubic8_coeffs})} and {(\ref{eq:u_coeff_repr})}. For this $u(x)$, at least one irreducible quartic polynomial is a factor of $\Phi_8(u(x))$.
\end{theorem}
  For a cubic polynomial $u(x)$ given by Theorem~\ref{thm:cubic8_lineq}, we take an irreducible quartic factor $r(x)$ from the factorization of $\Phi_8(u(x))$. If we let $t(x) \equiv u(x)^{2n+1}+1 \pmod{r(x)}\,\,(0 \le n \le 3)$, then Step~\ref{extbw:selt} of Algorithm~\ref{algo:Ext_BW} are finished.
  We continue the computation under the assumption that $k=8$ and $\sqrt{-D} \in \Q(\zeta_k)$. We can choose the CM discriminant $D=1$. Then we take $s(x)=(t(x)-1)^2 \mapsto \sqrt{-1}$ and execute Steps~\ref{extbw:construct}, \ref{extbw:check} in Algorithm~\ref{algo:Ext_BW} to get a family of curves.

  We now state the refinement of Algorithm~\ref{algo:Ext_BW} with restricted to our special case:
\begin{algorithm}
\label{algo:cubic8} Let $k=8, D=1, \deg u(x)=3$.  
\begin{description}
  \item[Input] A finite subset $S \subset \Q(\zeta_8)$.
  \item[Output] Families of elliptic curves with parameters $t(x),r(x),q(x)$.
\end{description}
  \begin{enumerate}[{\rm Step 1.}]
    \item \label{res:selomega} For each $\omega \in S$, compute $d, n_i$ by the equation~(\ref{eq:cubic8_invariant}), and let $S'=\{\omega \in S \mid d \not= 0, n_3 \not=0 \}$. If $S'$ is an empty set, then the algorithm fails.
    \item \label{res:selu} For each $\omega \in S'$, compute $u(x)$ by the equations (\ref{eq:cubic8_coeffs}) and (\ref{eq:u_coeff_repr}).
    \item \label{res:factor} For each $u(x)$ of Step \ref{res:selu}, compute all irreducible factors $r(x)$ of the polynomial $\Phi_8(u(x))$.
    \item \label{res:selt} For each pair $u(x),r(x)$ of Step \ref{res:factor}, compute all polynomials $t(x) \in \Q[x]$ such that $\deg t(x)<\deg r(x)$ and $t(x) \equiv u(x)^m+1 \pmod{r(x)}$ for all $m=1,3,5,7$.
    \item \label{res:sely} Compute $y(x) \equiv (2-t(x))(t(x)-1)^2 \pmod{r(x)} \quad (\deg y(x)<\deg r(x))$.
    \item \label{res:construct} Let $q(x)=(t(x)^2+y(x)^2)/4$.
    \item \label{res:check} For each triple $r(x),t(x),q(x)$ at Step~\ref{res:construct}, check whether $q(x),r(x)$ represent primes. If $q(x),r(x)$ represent primes, output a family $t(x),r(x),q(x)$.
  \end{enumerate}
\end{algorithm}

\subsection{Examples}
\label{sec:Res_FamilyExample}
\begin{table}[!h]
\caption{Sporadic families generate from cubic $u(x)$ with embedding degree 8}
\label{tab:family_cubicdeg8}
\begin{center}
  \begin{tabular}{|c|c|c|c|c|c|}\hline
    $\lc(u)$ & $u(x)$ & $t(x)$ & $\deg r(x)$ & $\deg q(x)$ & $\rho(t,r,q)$ \\ \hline \hline
    {\bf 2} & $\mathbf {2x^3+4x^2+6x+3}$ &$\mathbf{u(x)^3+1}$ & $\mathbf{4}$ & $\mathbf{6}$ & $\mathbf {3/2}$ \\ \hline
    9 & $9x^3+3x^2+2x+1$ &$u(x)^5+1$& $4$ & $6$ & $3/2$ \\ \hline
    17 & $17x^3+32x^2+24x+6$ & $u(x)^3+1$ & $4$ & $6$ & $3/2$ \\ \hline
    18 & $18x^3+39x^2+31x+7$ & $u(x)^3+1$ & $4$ & $6$ & $3/2$ \\ \hline
    64 & $64x^3+112x^2+75x+18$& $u(x)^5+1$ & $8$ & $14$ & $7/4$ \\ \hline
    68 & $68x^3+110x^2+65x+15$& $u(x)^5+1$ & $4$ & $6$ & $3/2$ \\ \hline
    {\bf 82} & $\mathbf {82x^3+108x^2+54x+9}$ &$\mathbf{u(x)^5+1}$ & $\mathbf{4}$ & $\mathbf{6}$ & $\mathbf {3/2}$ \\ \hline
    144 & $144x^3+480x^2+539x+202$&$u(x)^5+1$& $8$ & $14$ & $7/4$ \\ \hline
    144 & $144x^3+96x^2+29x+2$ &$u(x)^5+1$& $8$ & $14$ & $7/4$ \\ \hline
    216 & $216x^3+372x^2+263x+69$& --- & --- & ($\dagger$) & --- \\ \hline
    225 & $225x^3+2x$ & ---& --- & ($\dagger$) & --- \\ \hline
    257 & $257x^3+256x^2+96x+12$ & $u(x)^3+1$ & $4$ & $6$ & $3/2$ \\ \hline
    388 & $388x^3+798x^2+561x+134$ & $u(x)^5+1$ & $4$ & $6$ & $3/2$ \\ \hline
    392 & $392x^3+980x^2+821x+231$ & $u(x)^5+1$ & $8$ & $14$ & $7/4$ \\ \hline
    450 & $450x^3+11x$ & --- & --- & ($\dagger$) & --- \\ \hline
    {\bf 626} & $\mathbf{626x^3+500x^2+150x+15}$ &$\mathbf{u(x)^5+1}$ & $\mathbf{4}$ & $\mathbf{6}$ & $\mathbf {3/2}$ \\ \hline
    {\bf 738} & $\mathbf{738x^3+1488x^2+1006x+229}$ &$\mathbf{u(x)^5+1}$ & $\mathbf{4}$ & $\mathbf{6}$ & $\mathbf {3/2}$ \\ \hline
    800 & $800x^3+9x$ &$u(x)^5+1$& $8$ & $14$ & $7/4$ \\ \hline
    873 & $873x^3+969x^2+379x+53$ & $u(x)^7+1$ & $4$ & $6$ & $3/2$ \\ \hline
  \end{tabular}
\end{center}
\end{table}

  In Table~\ref{tab:family_cubicdeg8} We give a result of computations of the polynomial $u(x)$ by MAGMA \cite{MAGMA} using Algorithm~\ref{algo:cubic8}. The heading $\lc(u)$ denotes the leading coefficient of $u(x)$. We choose the input $S=\{\omega \in \Q(\zeta_8) \mid \omega=\sum_{i=0}^{3}a_i x^i ,\,a_i \in \Z,\, 0 \le a_i \le 300 \}$. 
%In order to make the polynomial coefficients small, we replaced $u(x)$ obtained in Step \ref{res:selu} of Algorithm~\ref{algo:cubic8} by $u(ax+b) \in \Z[x]$ with certain $a,b \in \Q, a\not=0$ . 
In the actual computation to make polynomial coefficients small, we further transform $u(x)$ obtained by Step \ref{res:selu} of Algorithm~\ref{algo:cubic8} to $u(ax+b) \in \Z[x]$ with suitable $a, b \in \Q, a\not=0$. 
%$ \in \Z[x] \,\, (a,b \in \Q, a\not=0)$ that has small coefficient instead of $u(x)$ given by Step  on . %$B$3$NA`:n$K$h$j(BStep \ref{res:factor}$B$K$*$$$F(B$r(x) \in \Z[x]$$B$G$+$D(B$\lc(r)$$B$,>.$5$/$J$k$3$H$,4|BT$5$l$k!%(B
  After computation by MAGMA, we tried to construct families for $\lg \lc(u)<10$. 
%$B%"%k%4%j%:%`(B~\ref{algo:cubic8}$B$K=>$$(B$(t,r,q)$$B$r7W;;$7!$Dj5A(B~\ref{def:pair_friendly_family}$B$rK~$?$7!$(B$\deg r(x)=4$$B$rK~$?$9AH$rM%@h$7$F:NMQ$7$?!%(B
  We explain the symbols in Table~\ref{tab:family_cubicdeg8}. For the column $\deg q(x)$, the symbol ($\dagger$) denotes that $q(x)$ does not {\it represent primes} for all pair $t(x),\,r(x)$. For rows, the {\bf bold notation} means that there exists a family of curves such that both $q(x)$ and $r(x)$ are primes for many integers $x$.

  We discovered many pairing-friendly families of curves with $\rho=3/2$ and also rediscovered a family which has $\lc(u)=9$ by Freeman et al. \cite[Example 6.18]{FST06}. It is interesting to note that for $\lc(u)=2,\, 82,\, 626,\, 738$, both $q(x)$ and $r(x)$ are primes for (infinitely) many integers $x$. We describe the case where $\lc(u)=82$ in detail as follows.
\begin{theorem}
\label{thm:cubic8_t82}
  The polynomials $t(x),\,r(x),\,q(x) \in \Z[x]$ given as follows represent a family of elliptic curves with embedding degree $k=8$ and the CM discriminant $D=1$. This family indeed generates pairing-friendly elliptic curves.
\[
  \begin{array}{rcl}
    t(x) & = & -82x^3-108x^2-54x-8 \\
    r(x) & = & 82x^4+108x^3+54x^2+12x+1 \\
    q(x) & = & 379906x^6+799008x^5+705346x^4\\
         &   &  +333614x^3+88945x^2+12636x+745
  \end{array}
\]
\end{theorem}
\begin{proof}
  The former half is already proved by Algorithm~\ref{algo:cubic8}, so we only need to prove the latter half. We may verify both $q(x_0)$ and $r(x_0)$ are primes with some integer $x_0$. We take $x_0=104$, then we get $q(x_0)=490506332802458249$ and $r(x_0)=9714910817$. Both of these are primes, so we can generate pairing-friendly curves by them.
\qed\end{proof}
  From a family of curves, we can actually construct pairing-friendly curves. Find an integer $x$ such that $q(x)$ is a prime and check whether $r(x)$ is a prime. To find such an integer $x$, we can reduce the number of the candidate by the chinese remainder theorem.
\begin{lemma}
\label{lem:cubic8_t82_4mod10}
  If an integer $q(x)$ in Theorem~\ref{thm:cubic8_t82} is a prime, then $x \equiv 14, 24\pmod{30}$.
\end{lemma}
\begin{proof}
  We can easily check that all $q(x)$ is even if $x$ is odd. We see that %look at constant of $q(x)$, then compute $q(x) \bmod 5$,
  \[
    q(x) \equiv x^6+3x^5+x^4+4x^3+x \pmod{5} .
  \]
  So $q(x) \equiv 0 \pmod{5}$ if $x \not\equiv 4 \pmod{5}$. In the same way we see that
  \[
    q(x) \equiv x^6+x^4+2x^3+x^2+1 \pmod{3} .
  \]
 So $q(x) \equiv 0 \pmod{5}$ if $x \not\equiv 1 \pmod{3}$. Then by the Chinese remainder theorem, $q(x)$ has no prime factor $2,3$ and $5$ only if $x \equiv 14, 24 \pmod{30}$. 
\qed\end{proof}

\section{Examples of pairing-friendly curves}
\label{sec:Res_Comp_Result}

  By Theorem~\ref{thm:cubic8_t82}, we can generate pairing-friendly curves using \cite[Theorems 3,4]{Mor98}. The elliptic curve $E/\Fq$ with the CM discriminant $D=1$ is represented as 
\[
  E: Y^2 \equiv X^3 + aX \pmod{q} \quad (a\not=0)
\]
where $a$ is parameter. Since $t$ is always divided by $4$ from the form of $t(x)$ in Theorem~\ref{thm:cubic8_t82}, we can easily compute $a$ by the method described in \cite{Mor98}. Using this, we give some numerical examples.
\begin{example}
   For $x=24000000000010394$ $(\lg x \approx 54.4)$, we get
  \[\begin{array}{rcl}
       q &  = & 726011672004446604951703464791789328991217313776602768811 \\
         &    & 50532069758156754787842298703647640196322590069, \\
       r &  = & 272056320000471307161600306182614014808404525177076771934 \\
         &    & 82845476817 \quad \mbox{(224-bit)},\\
       t &  = & -1133568000001472850432000637893917136092090964291460, \\
  \#E(\Fq)&  = & 726011672004446604951703464791789328991217313776602780147 \\
         &    & 18532071231007186788480192620783732287286881530, \\
       a &  = & 363005836002223302475851732395894664495608656888301384405 \\
         &    & 75266034879078377393921149351823820098161295035  .\\
    \end{array}\]
  Then $\lg r \approx 224.0$, $\lg q \approx 345.0$ and $\rho(E) \approx 1.54$.
\end{example}

\begin{example}
  For $x=6130400000000029634$ $(\lg x \approx 62.4)$, we get
  \[
    \begin{array}{rcl}
       q &  = & 20165501539097468598089799012338448337497685\\
         &    & 26807341931299469596014851929961512795928195\\
         &    & 2496431544631024161702159356789, \\
       r &  = & 11581614432149089047832789189966585476390503\\
         &    & 3269185946585920376349372307631217 \quad \mbox{(256-bit)},\\
       t &  = & -1889210236224232197405821630084439441516429\\
         &    & 1734047019380020, \\
  \#E(\Fq)&  = & 20165501539097468598089799012338448337497685\\
         &    & 26807341931299471485225088154193710201749825\\
         &    & 3340825959795315895749178736810, \\
       a &  = & 10082750769548734299044899506169224168748842\\
         &    & 63403670965649734798007425964980756397964097\\
         &    & 62482157723155120808510796783952  .\\
    \end{array}
  \]
  Then $\lg r \approx 256.0$, $\lg q \approx 393.0$ and $\rho(E) \approx 1.54$.
\end{example}

\begin{example}
   For $x=-72057594037930756$ $(\lg x \approx 56.0)$, we get
  \[
    \begin{array}{rcl}
       q &  = & 5318077912637504134292767901251647400395578540 \\
         &    & 3827730100050941212371435046023372666628598916 \\
         &    & 049952969199369, \\
       r &  = & 2210715626706698491377041180063927762099958931 \\
         &    & 722603805474805907424817 \quad \mbox{(230-bit)}, \\
       t &  = & 3067984237085391549834039420816298507616442947 \\
         &    & 7994640,\\
  \#E(\Fq)&  = & 5318077912637504134292767901251647400395578540 \\
         &    & 3827730069371098841517519547682978458465613839 \\
         &    & 885523491204730, \\
       a &  = & 1772692637545834711430922633750549133465192846 \\
         &    & 7942576700016980404123811682007790888876199638 \\
         &    & 683317656399790  .\\
    \end{array}
  \]
  Then $\lg r \approx 230.4$, $\lg q \approx 354.5$ and $\rho(E) \approx 1.54$. 
For the Ate pairing \cite{HSV06}, it is important that $t$ has a low hamming weight for computation. We tried to find a curve with $r$ between 224 bit and 256 bit, we found that $r$ has a Hamming weight 72 and $t$ has a Hamming weight 45 in this example.
\end{example}

\section{Conclusion}
\label{sec:Cons}

  We proposed a new algorithm for systematically constructing families of elliptic curves with given embedding degree and the CM discriminant. It was shown to be efficient by producing actual families of curves and explicit numerical examples for the case of embedding degree 8. The key point is employing the method of indeterminate coefficients to choose polynomials. Obviously our method of indeterminate coefficients are also applicable to the general case. % Our numerical examples are not the curves of prime order different from the cases of \cite{BN06, Fre06}, but they are %We only give a one of the most realistic solution, but not best, such as Barreto-Naehrig family with $k=12$  or Freeman curves with $k=10$ \cite{Fre06,FST06}. We will apply the results of a study which found to hyperelliptic curves, or experiment in new cryptographic schemes. Of course, we continue to research best construction of family of curves.

\bibliographystyle{splncsplain}
\bibliography{ecpf}

\begin{thebibliography}{10}

\bibitem{AM93}
Atkin, A.O.L., Morain, F.:
\newblock Elliptic curve and primarity proving.
\newblock Math. Comp. \textbf{61}(203) (1993)  29--68.

\bibitem{BK98}
Balasubramanian, R., Koblitz, N.:
\newblock The improbability that an elliptic curve has subexponential discrete
  log problem under the menezes-okamoto-vanstone algorithm.
\newblock Journal of Cryptology \textbf{11}(2) (1998)  141--145.

\bibitem{BN06}
Barreto, P.S.L.M., Naehrig, M.:
\newblock Pairing-friendly elliptic curves of prime order.
\newblock In: SAC2005 - Workshop on Selected Areas in Cryptography. Volume 3897
  of Lecture Notes in Computer Science., Springer (2006)  319--331.

\bibitem{BF03}
Boneh, D., Franklin, M.:
\newblock Identity based encryption from the weil pairing.
\newblock SIAM Journal of Computing \textbf{32}(3) (2003)  586--615.

\bibitem{BLS01}
Boneh, D., Lynn, B., Shacham, H.:
\newblock Short signatures from the weil pairing.
\newblock In: Advances in Cryptology - AsiaCRYPT 2001. Volume 2248 of Lecture
  Notes in Computer Science., Springer (2001)  514--532.

\bibitem{BW05}
Brezing, F., Weng, A.:
\newblock Elliptic curves suitable for pairing based cryptography.
\newblock Designs, Code and Cryptography \textbf{37}(1) (2005)  133--141.

\bibitem{FST06}
Freeman, D., Scott, M., Teske, E.:
\newblock A taxonomy of pairing-friendly elliptic curves.
\newblock preprint (2006).
\newblock http://math.berkeley.edu/~dfreeman/papers/taxonomy.pdf.

\bibitem{GMV05}
Galbraith, S., Mckee, J., Valen\c{c}a, P.:
\newblock Ordinary abelian varieties having small embedding degree.
\newblock In: Workshop on Mathematical Problems and Techniques in Cryptology,
  Barcelona, CRM (2005)  29--45.

\bibitem{HSV06}
Hess, F., Smart, N., Vercauteren, F.:
\newblock The eta pairing revisited.
\newblock Cryptology ePrint Archive, Report 2006/110 (2006).
\newblock http://eprint.iacr.org/2006/110/.

\bibitem{Jou04}
Joux, A.:
\newblock A one round protocol for tripartite diffie-hellman.
\newblock Journal of Cryptology \textbf{17}(4) (2004)  263--276.

\bibitem{MAGMA}
{MAGMA Group}:
\newblock Magma computational algebra system.
\newblock http://magma.maths.usyd.edu.au.

\bibitem{Mor98}
Morain, F.:
\newblock Primality proving using elliptic curves: An update.
\newblock In: ANTS 1998 - 3rd Algorithmic Number Theory Symposium. Volume 1423
  of Lecture Notes in Computer Science., Springer (1998)  111--127.

\bibitem{SOK00}
Sakai, R., Ohgishi, K., Kasahara, M.:
\newblock Cryptosystems based on pairing.
\newblock In: 2000 Symposium on Cryptography and Information Security
  (SCIS2000). (2000).

\end{thebibliography}

\end{document}